\documentclass[10pt,twocolumn]{IEEEtran}

\usepackage[latin9]{inputenc}
\usepackage{verbatim}
\usepackage{float}
\usepackage{graphicx,epsfig,psfrag}
\usepackage{url,color}
\usepackage[cmex10]{amsmath}
\usepackage{subfigure,graphics}
\usepackage{mdwmath,amsfonts}
\usepackage[british]{babel}
\usepackage[colorinlistoftodos,bordercolor=orange,backgroundcolor=orange!20,linecolor=orange,textsize=scriptsize]{todonotes}
\usepackage{enumitem}

\newtheorem{theorem}{Theorem}

\newtheorem{proposition}{Proposition}
\newtheorem{lemma}{Lemma}
\newtheorem{definition}{Definition}

\newtheorem{remark}{Remark}

\newtheorem{prob}{Problem}
\newtheorem{alg}{Algorithm}

\usecounter{algorithmenumi}

\makeatletter

\IEEEoverridecommandlockouts

\DeclareMathOperator{\diag}{diag}

\DeclareMathOperator{\opd}{d\!}

\newcommand{\A}{\mathcal{A}}
\newcommand{\B}{\mathcal{B}}
\newcommand{\C}{\mathcal{C}}
\newcommand{\K}{\mathcal{K}_e}
\newcommand{\W}{\mathcal{W}}
\newcommand{\M}{\mathcal{M}}

\newcommand{\Gr}{\mathcal{G}}
\providecommand{\rset}[1]{\mathbb{R}^}

\makeindex

\title{$H_2$ model reduction of linear network systems by moment matching and optimization}

\author{I. Necoara~\thanks{I. Necoara is with the Department of Automatic Control and Systems Engineering,  University Politehnica Bucharest, 060042 Bucharest, Romania. E-mail: \tt{ion.necoara@acse.pub.ro}.} and T.C. Ionescu~\thanks{T.C. Ionescu is with the Department of Automatic Control and Systems Engineering, University Politehnica Bucharest, 060042 Bucharest, Romania \&  Institute of Mathematical Statistics and Applied Mathematics of the Romanian Academy, 050711 Bucharest, Romania. E-mail: {\tt tudor.ionescu@acse.pub.ro}.} 
\thanks{This work  is supported by the Executive Agency for Higher Education, Research and Innovation Funding (UEFISCDI), Romania,  PNIII-P4-PCE-2016-0731, project ScaleFreeNet, no. 39/2017. } }

\makeatother

\begin{document}
\markboth{IEEE Trans. Autom. Control, December 2018}{}
\maketitle

\begin{abstract}
In this paper we study the problem of model reduction of linear network systems. We aim at computing a reduced order stable approximation of the network with the same topology  and optimal w.r.t. $H_2$-norm error approximation. Our approach is based on time-domain moment matching framework, where we optimize over families of parameterized reduced order models  matching a set of moments at  arbitrary interpolation points. The  parameterization of the low order models is in terms of the free parameters and of  the interpolation points.  For this family of parameterized models  we formulate an optimization-based  model reduction problem with  the $H_2$-norm of error approximation as  objective function  while the preservation of some structural and physical properties yields  the constraints. This problem is nonconvex and we write it in terms of the Gramians of a minimal realization of the error system. We propose  two solutions for this problem. The first solution assumes that the error system  admits a block diagonal observability Gramian, allowing for a simple convex  reformulation as semidefinite programming, but  at the cost of some performance loss.  We also derive sufficient conditions to guarantee block diagonalization of the Gramian. The second solution employs a gradient projection method for a smooth reformulation yielding   (locally) optimal  interpolation points and  free parameters. The potential of the methods is illustrated on several network examples.
\end{abstract}

\setlength{\arraycolsep}{2pt}

\section{Introduction}
\label{intro}

\IEEEPARstart{C}{omplex network systems} consist of multiple interacting dynamical subsystems, interconnected through a graph enabling the subsystems to share information,  coordinate their activities and have self-control mechanisms \cite{NecNed:11}. However, the
corresponding models of  network systems are too complex and difficult to analyze, rendering it is almost impossible to systematically develop operating and/or open/closed-loop control algorithms. Therefore, we need approximation models to do analysis, simulation and control. 

\medskip

\noindent \textit{State-of-the-art}: The problem of model  reduction of interconnected systems has  been long studied in different frameworks, see e.g. \cite{reis-stykel-SPRINGER2008} and references therein for a survey. There are two main existing approaches.  

\noindent A {\em first} approach stemming from mathematics considers network systems as static mathematical objects. The reduction is treated with topological objectives, focusing on obtaining a reduced network abstracting a large-scale network by merging groups of nodes into super-nodes (so-called clustering) \cite{besselink-sandberg-johansson-ECC2014}.  For example, \cite{MonTre:13} aims at preserving stability and synchronisation of the system, \cite{CheKaw:16a} preserves an interconnection structure and  synchronization   by  aggregating subsystems with similar frequency responses, \cite{IshKas:14} provides a reduced system with a dynamical behaviour close to the initial system while preserving several properties for control purpose, and \cite{MarFra:18} proposes a  network reduction method  preserving the   flow network property and  the reduced  graph to be scaled-free.  

\noindent A {\em second} approach comes from systems and control theory. The aim is to reduce the network system by preserving  consistency/structure in the network. Here, one category of results are in the framework of stability preserving  balanced truncation, see e.g. \cite{sandberg-murray-OCAM2009,sandberg-TAC2010}, where the balancing  yields the so-called structured Hankel singular values (invariants showing the importance of subsystems states with respect to a chosen input-output map for the whole  network). The states with the lowest structured Hankel singular values are truncated directly from the full model, resulting in a low order stable network  satisfying the given interconnection map.  A second category of results are based on interpolatory methods \cite{bai-ANM2002}, as  e.g. in  \cite{lutowska-2012,reis-stykel-MCDMS2007,vandendorpe-vandooren-2009}, see also \cite{li-bai-CMS2005,freund-ACM2004} for earlier results. Here, structured Krylov projections are  applied directly on the  entire network to preserve the topology. Note that none of the presented results reduce the number of nodes (subsystems) or alter the interconnection map of the network.

\medskip

\noindent \textit{Motivation}:   
In this paper we also consider the second approach of approximating the subsystems of a network. If the number of subsystems is large, we first reduce their number using clustering techniques  and then perform subsystem approximation. To the best of our knowledge, in the time-domain moment matching framework \cite{astolfi-TAC2010,i-astolfi-colaneri-SCL2014}, finding a reduced model  of the network optimally w.r.t. the $H_2$-norm of the  approximation error  in the family of $\nu$ order models that matches a set of $\nu$  moments, while preserving the network topology and stability, is an open question.  Some initial  progress has been made recently in \cite{NecIon:18} for general linear  systems. However, a direct application of this approach on the full model of the network system does not preserve network topology. Hence, this unsolved problem motivates our work here. 

\medskip

\noindent \textit{Contributions}:    In this paper we provide a systematic procedure for  approximating the subsystems of  a network  optimally  while preserving the network topology and  stability. The proposed procedure is based on time-domain moment matching, where families of parametrized low order models matching a set of moments at arbitrary interpolation points are computed. Here we use  the free parameters and the interpolation points defining the parameterization to find the optimal approximation of the network measured in terms of the   $H_2$-norm of the error system.  We formulate an optimization problem with the $H_2$-norm of the approximation error as objective function, while the preservation of some structural and physical properties yields  the constraints. The problem is nonconvex and we prove that it can be written  in terms of the controllability/observability Gramians of a minimal realization of the error system. We propose  two solutions for solving this problem. The first solution assumes that the error system  admits a block diagonal observability Gramian, allowing for a simple convex  reformulation as semidefinite programming, at the cost of some performance loss.  Also,  sufficient conditions are derived to guarantee block diagonal Gramians.  The second solution employs a gradient projection method for a smooth reformulation  yielding   $H_2$ (locally)  optimal  reduced order  models.   Both solutions provide (optimal) stable low order network models,  parameterized in the interpolation points and in the free parameters, matching a  set of  moments and preserving the interconnection map of the original network. The efficiency of the methods is illustrated on a  positive network  and on a multiple area power system.

\medskip

\noindent \textit{Content}:  The paper is organized as follows. In Section \ref{sect_prel} we briefly review the time-domain moment matching model reduction of linear  systems. In Section \ref{sect_optimalMOR} we formulate the optimization-based $H_2$-norm moment matching model reduction problem with a Gramian-type cost function. In Sections \ref{sect_SDP} and \ref{sect_gradient}, we propose two  numerical optimization methods for solving this model reduction problem and extensions are given in Section \ref{sect_ext}. In Section \ref{sect_exmp} we illustrate the efficiency of our theory on several network examples.

\medskip

\noindent \textit{Notation}:   $\mathbb{R}$ and $\mathbb{C}$ denotes the set of real and complex numbers, respectively.  For a positive integer $N$ we denote by $[N] =\{1, \cdots, N\}$. For a matrix $A \in \mathbb{R}^{n\times n}$,  $\sigma(A)$ denotes the set of its eigenvalues and $A_{ij}$ indicates the $(i,j)$  matrix block of $A$ of appropriate dimension.


\section{Preliminaries}
\label{sect_prel}
\noindent In this section,  we briefly present  the main results on time-domain moment
matching for  linear  systems \cite{i-astolfi-colaneri-SCL2014}.

\subsection{Moments and moment matching
\label{subsec:Moment_def}}
\noindent Consider the linear time-invariant system:
\begin{align}
\label{system}
\dot x & = Ax + Bu \quad \text{and} \quad
y = Cx,
\end{align}
where $x \in \mathbb{R}^n$ is the state of the system,
$u \in \mathbb{R}^m$ is the input and $y \in \mathbb{R}^p$
is the output, respectively. Consequently, system matrices $A \in \mathbb{R}^{n \times n}, B \in  \mathbb{R}^{n \times m}$ and $C \in \mathbb{R}^{p \times n}$.  Throughout the paper we assume that the system is stable (i.e. $\sigma(A) \subset \mathbb C^-$) and that \eqref{system} is a minimal (i.e. controllable and observable) realization of the transfer function:
\begin{equation}
\label{tf}
K(s) =C(sI - A)^{-1} B.
\end{equation}

\noindent The moments of linear system \eqref{system} at a point $s \in \mathbb{C}$ on the complex plane are defined as follows:
\begin{definition}
\label{def_moment}\cite{antoulas-2005,astolfi-TAC2010,gallivan-vandendorpe-vandooren-SIAM2004} The $k$-moment of the system (\ref{system}) at the point $s\notin\sigma(A)$, along the direction $\ell\in\mathbb C^m$ is  $\eta_{k}(s;\ell)={(-1)^{k}}/{k!} \ {\opd^k \! K(s)}/{\opd s^{k}} \ell \in \mathbb C^{p},$ with $k\geq 0$ integer.
\end{definition}

\noindent Consider the matrix $S \in \mathbb{R}^{\nu\times\nu}$, with $\nu \leq n$ and  $\sigma(S)=\{s_i: i\geq 0 \; \text{integer}\}  \subset \mathbb C$. Let $L =[\ell_1 \cdots   \ell_\nu] \in\mathbb{R}^{m \times \nu}$ be such that the pair $(L,S)$ is observable. Since the system is assumed minimal, the Sylvester equation:
\begin{equation}
A\Pi+BL=\Pi S,\label{eq_Sylvester_Pi}
\end{equation}
has the unique solution $\Pi\in \mathbb{R}^{n\times\nu}$ with ${\rm rank}\ \Pi=\nu$  provided that $\sigma(A)\cap\sigma(S)=\emptyset$ \cite{antoulas-2005}.  Then, the moments of a system can be characterized as follows:

\begin{proposition} \cite{i-astolfi-colaneri-SCL2014}
\label{def_PI} Consider the system \eqref{system} and let $\Pi$ be the unique solution of equation (\ref{eq_Sylvester_Pi}). Then, at the interpolation points $s_i\in\sigma(S)$, the moments of the system (\ref{system}) along directions $\ell_i$, $\eta_{k}(s_{i};\ell_i)$, with $i,k\geq 0$ integers,  are characterized by the  matrix~$C\Pi$.
\end{proposition}

\noindent We now  present the moment matching property and  the reduced order model satisfying it:
\begin{proposition}
\cite{astolfi-TAC2010,i-astolfi-colaneri-SCL2014}
\label{thm_MMcond}
Consider the $\nu$ order linear system:
\begin{align}
\label{red_mod_F}
\dot \xi & = F \xi + Gu \quad \text{and} \quad
\psi = H \xi,
\end{align}
with the state $\xi \in \mathbb{R}^\nu$, input $u \in \mathbb{R}^m$  and output $\psi \in \mathbb{R}^p$.   Here,  $\nu \leq n$,  $F \in \mathbb{R}^{\nu \times \nu}, G \in  \mathbb{R}^{\nu \times m}$ and $H \in \mathbb{R}^{p \times \nu}$.  Assuming $\sigma(F)\cap\sigma(S)=\emptyset$,  then the reduced order system \eqref{red_mod_F} matches the moments of \eqref{system} at $\sigma(S)$ if and
only if:
\begin{equation}\label{cond_MM}
HP=C\Pi,
\end{equation}
where the invertible matrix $P \in \mathbb{R}^{\nu \times \nu}$ is
the unique solution of the Sylvester equation
$FP+GL=PS.$
\end{proposition}

\noindent Note that the invertible matrix $P$ in Proposition \ref{thm_MMcond} is merely a
coordinate transformation. Hence, taking $P=I_\nu$ yields a  parameterized $\nu$ order model (with the free parameters $(G,L)$ and the interpolation points  matrix $S$) achieving moment matching at $\sigma(S)$, as shown in the next result:
\begin{proposition}
\cite{astolfi-TAC2010,i-astolfi-colaneri-SCL2014}
\label{thm_redmod_CPi}
Assume that $(L,S)$ is observable and  $\sigma(A)\cap\sigma(S)=\emptyset$.
Consider  the $\nu$ order linear system:
\begin{align}
\label{redmod_CPi}
\widehat \Sigma_{(S,G,L)}: \quad\dot \xi & = (S- GL) \xi + Gu, \quad
\psi = C \Pi \xi,
\end{align}
with  $\nu \leq n$ and  the transfer function:
\begin{equation}\label{tf_redmod_CPi}
\widehat K(s)=C\Pi(sI-S+GL)^{-1}G,
\end{equation}
where $\Pi$ is the unique solution of (\ref{eq_Sylvester_Pi}). Assuming that
$\sigma(S-GL)\cap\sigma(S)=\emptyset$,  then the system (\ref{redmod_CPi}), with the transfer function \eqref{tf_redmod_CPi},  is a reduced order model of (\ref{system}) parametrized in  $S, G$ and $L$,  matching the moments $C \Pi$ of  system \eqref{system} at $\sigma(S)$.
\end{proposition}

\begin{remark}
\label{rem_Lfix}
The system $\widehat \Sigma_{(S,G,L)}$ in \eqref{tf_redmod_CPi} describes a  $\nu$ order approximation of \eqref{system}  that \emph{achieves moment matching} at $\sigma(S)$.   Since $L$ is only used in  $\widehat \Sigma_{(S,G,L)}$  to ensure observability of the pair $(L,S)$ and since observability is  generic,  then,  without loss of generality,  we  fix \emph{a priori} matrix $L$, i.e., we  fix the directions $\ell_i$ to compute moments along, see e.g. \cite{astolfi-TAC2010, i-astolfi-colaneri-SCL2014}.  Hence, in the rest of the paper we consider $\widehat\Sigma_{(S,G,L)}=\widehat \Sigma_{(S,G)}$, defining a \textit{family} of $\nu$ order models  matching $\nu$ moments along  \textit{fixed} directions $\ell_i$ of system   \eqref{system} at $\sigma(S)$, for all $G$,  such that:
\begin{enumerate}[wide,nosep,label=\roman*)]
\item $\widehat\Sigma_{(S,G)}$ is parametrized in $(S,G)$
\item $\sigma(S-GL) \cap \sigma(S)=\emptyset$.
\end{enumerate}	
\end{remark}


\subsection{$H_2$-norm based on the Gramians of linear systems}
\noindent We now  briefly recall the definition and computation of the $H_2$-norm of a linear system.  For   {\em minimal stable} system \eqref{system} with  transfer function  \eqref{tf},  the ${H}_2$-norm is defined as \cite{gugercin-antoulas-beattie-SIAM2008}:
\[  \|K\|_{{H}_2}=\sqrt{\int_{-\infty}^{\infty} |K(j\omega)|^2 \opd\omega}.  \]
This norm can be written explicitly in matrix form as \cite{gugercin-antoulas-beattie-SIAM2008}:
\begin{equation}
\label{eq_2norm_gram}
\|K\|_{{H}_2}^2=C^TWC=B^TMB,
\end{equation}
where $W$ is the controlability Gramian and $M$ is the observability Gramian of the linear system  \eqref{system}.


\section{Optimal $H_2$ model reduction formulation of linear network systems}  \label{sect_optimalMOR}
\noindent In this section we formulate a model reduction problem, yileding a family of parametrized models for each subsystem of a linear network without altering its structure. To determine the best approximation in terms of the $H_2$-norm of the error, we propose an optimization formulation, with the $H_2$-norm of the error as objective function, while stability and structure are imposed as constraints.

\subsection{Linear network systems}
\noindent We  perform model reduction for linear network systems consisting of $N$ interconnected subsystems, with dynamics defined by the  linear state
space equations:
\begin{align}
\label{mod3} 
\dot x_i & =  \sum _{j \in {\mathcal N}_i}  A_{ij} x_{j} + B_{i} u   \quad  \forall i \in [N],
\end{align}
where $x_i \in \mathbb{R}^{n_i}$  represents the state of the  $i$th subsystem,  $u \in \mathbb{R}^m$ is the common input, $A_{ii} \in \mathbb{R}^{n_i
\times n_i}$, $B_i \in \mathbb{R}^{n_i \times m}$, and $A_{ij} \in \mathbb{R}^{n_i \times n_j}$. The index set ${\mathcal N}_i \subseteq [N]$ contains the index $i$ and all
the indices of the subsystems which interact with the subsystem $i$. Thus, in  \eqref{mod3} we consider that each subsystem is influenced through the states of the neighboring subsystems. For a more general network description see  Section \ref{sect_ext}.

\begin{figure}[h!]
\centering{ \psfrag{S1}{$\Sigma_1$}   \psfrag{S2}{$\Sigma_2$}
\psfrag{S3}{$\Sigma_3$}    \psfrag{S4}{$\Sigma_4$}
\includegraphics[width=0.3\textwidth]{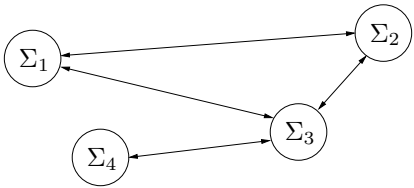} } 
\caption{An example of a network system composed of four subsystems.} 
\label{fig:net_systems}
\end{figure}
\noindent For example, consider  the network system in Figure \ref{fig:net_systems}, where the arrows indicate the interactions between the subsystems $\Sigma_1,\Sigma_2,\Sigma_3$ and $\Sigma_4$. If we consider the fourth subsystem $\Sigma_4$, we have ${\mathcal
N}_4=\{3,4\}$ and hence:
\[ \dot x_4  = A_{44}x_4 + A_{43} x_3 + B_4 u \quad \text{and} \quad A_{41} = A_{42} =0. \]
For model reduction, we also express the dynamics of the entire  network system in the compact form \eqref{system}, $\dot x  = Ax + Bu$, where  $x = [x_1^T \cdots x_N^T]^T \in \mathbb{R}^n$, with $n = \sum_{i=1}^N n_i$,  denotes the states of the entire network and the input  $u \in \mathbb{R}^m$.  As output of the  network system we consider a linear combination of the states of each subsystem:    $y = C x$,  \text{where}   $ C \in \mathbb{R}^{p \times n}$.    
Then, the system matrices  $A \in \mathbb{R}^{n \times n}$, $B \in
\mathbb{R}^{n \times m}$ and $C = [C_1 \dots C_N] \in \mathbb{R}^{p \times n}$, with $C_i  \in \mathbb{R}^{p \times n_i}$, are given by:
\begin{equation}
\left[\begin{array}{c|c} A & B \\
\hline C & 0
\end{array}\right] =\left[\begin{array}{cccc|c}
A_{11} & A_{12} & \dots & A_{1N} & B_{1}\\
A_{21} & A_{22} & \dots & A_{2N} & B_{2}\\
\vdots & \vdots & \ddots & \vdots & \vdots\\
A_{N1} & A_{N2} & \dots & A_{NN} & B_{N}\\
\hline C_{1} & C_{2} & \dots & C_{N} & 0
\end{array}\right],
\label{realization_block}
\end{equation}
where we have that the matrix block $(i,j)$ of $A$ satisfies:
\[  A_{ij} = 0 \quad  \forall i \in [N], \;  j \not \in {\mathcal N}_i. \]
Recall that we assume stable networks, i.e.  $\sigma(A) \subset \mathbb C^-$.  Note that the dimension $n$ of entire network system is usually too
large, so that it is almost impossible to develop  open or
closed-loop control algorithms in a systematic way. Therefore, we need to obtain approximation models of the network system \eqref{realization_block}, useful for analysis,
simulation and control. Unfortunately, moment matching-based model reduction techniques, such as in \cite{NecIon:18}, do not preserve the network structure. However, working with an approximation violating basic network constraints it always leaves the question of how conclusive the results on this basis~are.


\subsection{Optimal $H_2$ moment matching-based model reduction problem preserving network structure}
\noindent In this section, we formulate  the component model reduction problem of the
linear network system \eqref{mod3} with the network structure given
in \eqref{realization_block}. Recall that  $A_{ij} = 0$ in
\eqref{realization_block} if $i \in [N], \; j \not \in {\mathcal N}_i$ and the
dimension of the entire network system is $n = \sum_{i=1}^{N}n_{i}$.
The goal is to perform model order reduction such that the network
structure is preserved, i.e. compute reduced order models for each
subsystem of the form:
\begin{align}
\label{mod4} \dot \xi_i & =  \sum _{j \in {\mathcal N}_i}
F_{ij} \xi_{j} + G_{i} u   \qquad  \forall i \in [N],
\end{align}
where $\xi_i \in \mathbb{R}^{\nu_i}$, with $\nu_{i} \leq n_{i}$,
represents the reduced state of the  $i$th subsystem,  $u \in
\mathbb{R}^m$ is the common input, $F_{ij} \in \mathbb{R}^{\nu_{i} \times
\nu_{j}}$ and  $G_{i} \in \mathbb{R}^{\nu_{i} \times m}$.  Moreover, we want to preserve  the network structure, that is $F_{ij} = 0$ if $i \in [N], \; j \not \in {\mathcal N}_i$. If the number of subsystems $N$ is large, we first  reduce their number using  existing  clustering techniques (see Section I)  and then perform the subsystem approximation procedure described below. Note that the
dimension of the whole reduced model  $\dot \xi =  
F \xi + G u $  is $ \nu = \sum_{i=1}^{N}
\nu_{i}$, where $\xi =[\xi_1^T \cdots \xi_N^T]^T \in
\mathbb{R}^{\nu}$ denotes the full state of the reduced model. We also define the output of the reduced network:
\[  \psi = H \xi, \]
where $H \in \mathbb{R}^{p \times \nu}$. The  matrices of the reduced network system \eqref{mod4} are written in a compact
form as:
\begin{equation}
\left[\begin{array}{c|c} F & G \\
\hline H & 0
\end{array}\right] =\left[\begin{array}{cccc|c}
F_{11} & F_{12} & \dots & F_{1N} & G_{1}\\
F_{21} & F_{22} & \dots & F_{2N} & G_{2}\\
\vdots & \vdots & \ddots & \vdots & \vdots\\
F_{N1} & F_{N2} & \dots & F_{NN} & G_{N}\\
\hline H_{1} & H_{2} & \dots & H_{N} & 0
\end{array}\right],
\label{red_mod_F_block}
\end{equation}
where   $H_{j} \!\in\! \mathbb{R}^{p \times \nu_{j}}$ and consider also $L=[L_1  \dots L_N]$, with blocks $L_i \in \mathbb{R}^{m \times \nu_i}$.  We want $F_{ij} = 0$
 for  $i \in [N], \; j \not \in {\mathcal N}_i$.  Based on the parametrizations of the reduced model \eqref{redmod_CPi} given in terms of the  interpolation points matrix $S$ and of the free parameters $G$ (i.e. $F=S-GL, G$ and $H=C \Pi$) and using  the $H_2$-norm of the approximation error as objective function and physical and structural restrictions as constraints, we derive below an   optimization problem  to  determine the minimal approximation.  More precisely,  the optimal $H_2$ model reduction problem by moment matching is formulated as:

\begin{prob}\label{prob_optH2_general}\normalfont
Given a linear network system \eqref{system}  with the subsystem matrices \eqref{realization_block} and the transfer function $K(s)$ as in \eqref{tf}  and the directions of moments $L$, find a reduced order linear network $\widehat\Sigma_{(S,G)}$ of the form  \eqref{redmod_CPi} with subsystem matrices  \eqref{red_mod_F_block}  and  the transfer function $\widehat K(s)$  as in \eqref{tf_redmod_CPi},  parametrized in the interpolations matrix $S$ and  the free parameters $G$, that matches $\nu$ moments of \eqref{system} at $\sigma(S)$ and satisfies the constraints:
\begin{enumerate}[label=(\roman*),nosep]
\item the $H_2$-norm of the error system  $\|K-\widehat K\|_2$ is minimal
\item the reduced model $\widehat K$ is stable (i.e. $\sigma(S-GL) \subset \mathbb C^-$)
\item the matrix $F=S-GL$ preserves the network topology of $A$ (i.e for all $i \in [N]$,  $(S-GL)_{ij}=0$ if $j\notin{\cal N}_i$)
\item $\sigma(S) \cap \sigma(A) = \emptyset$,  $\sigma(S) \cap \sigma(S-GL) = \emptyset$ and the pair $(L,S)$ observable.
\end{enumerate}
\end{prob}

\noindent However,  it is difficult to deal with the restrictions (iv): $(S,L)$ observable,  $\sigma(S) \cap \sigma(A) = \emptyset$, and $\sigma(S) \cap \sigma(S-GL) = \emptyset$. One possibility  is to fix $S$ and $L$ such that  (iv) is automatically satisfied (e.g., without loss of generality, take $S = \diag(s_1,\dots,s_\nu)$, with $s_i \in \mathbb C^+ \; \forall i$, and  $L=[\ell_1  \dots \ell_\nu]$, with $\ell_i \not=0 \; \forall i$) and search only for the free parameters $G$.  All our results  hold  for this choice. Another possibility, which we also follow in this paper,  is to  fix $L$.  Note that since  model reduction procedures usually render  $S$  unstable, while $A$ and $S-GL$ are stable, the first two constraints in (iv) are automatically satisfied.  Moreover, since observability is generic, by Remark \ref{rem_Lfix}, also observability of $(L,S)$  holds. Hence, the constraints (iv) are not imposed in the numerical algorithms, but will be checked after yielding a solution to Problem \ref{prob_optH2_general}. 
%
%
Therefore, in the sequel we propose an optimization formulation of Problem \ref{prob_optH2_general}, without constraints (iv)   in   the unknowns  $S$ and $G$, while $L$ is fixed a priori.  Under these settings  Problem \ref{prob_optH2_general} can be recast  in terms of the Gramians of the realization of the \emph{error system}:
$$\K = K - \widehat K,$$ with $\widehat K$ from \eqref{tf_redmod_CPi}, parameterized in $(S,G)$. Let $(\A_e,\B_e,\C_e)$ be a state-space realization of the error transfer function $\K$: $$\K(s)=\C_e(sI-\A_e)^{-1}\B_e,$$ where
\begin{equation}\label{eq_err_realization}
\A_e=\begin{bmatrix} A & 0 \\ 0 & S - GL \end{bmatrix},\ \B_e=\begin{bmatrix} B \\ G \end{bmatrix},\ \C_e=C\begin{bmatrix} I &  -\Pi\end{bmatrix}.
\end{equation}

\noindent Denote the controllability and the observability Gramians of \eqref{eq_err_realization} by $\W$ and $\M$, respectively. They are solutions of the Lyapunov equations \cite{antoulas-2005}:
\begin{subequations}
\label{eq_Gramians_err}
\begin{align}
\label{W_Lyap}
\A_e\W+\W\A_e^T+\B_e\B_e^T &=0, \\
\label{M_Lyap}
\A_e^T\M+\M\A_e+\C_e^T\C_e &=0.
\end{align}
\end{subequations}
Since we assume $\sigma(A), \sigma(S-GL) \subset \mathbb C^-$, then  the matrix $\A_e$ is also stable, i.e., $\sigma(\A_e) \subset \mathbb C^-$. Hence,  there exist unique positive semidefinte solutions $\W$ and $\M$   of equations \eqref{eq_Gramians_err}, respectively. We partition  $\W$ and $\M$ following the two block structure of the error  matrix $\A_e$:
\begin{equation}
\W=\begin{bmatrix} W_{11} & W_{12} \\ W_{12}^T & W_{22}\end{bmatrix},\;\;   \M=\begin{bmatrix} M_{11} & M_{12} \\ M_{12}^T & M_{22}\end{bmatrix}.
\end{equation}

\noindent The communication graph between subsystems  imposes a constraint on the admissible parameterizations. If for subsystem $i \in [N]$,  $j \notin {\cal N}_i $, then there is no communication link from subsystem $j$ to subsystem $i$.
Thus, in the reduced model, the subsystem $i$ cannot be also influenced by the subsystem
$j$. This leads to the structured constraint $\Gr$,  defined by:
\[  \Gr= \{ (S,G): \;  (S-GL)_{ij} =0 \; \forall i \in [N], j \notin {\cal N}_i   \}.  \]
Let us define the feasible set for the reduced model:
\begin{align*}
{\cal R}= \left\{ (S, G): \;  \sigma(S-GL) \subset \mathbb C^-, \; (S,G) \in \Gr \right\}.
\end{align*}
By \eqref{eq_2norm_gram} and the simplification stated above,  Problem  \ref{prob_optH2_general} becomes:
\begin{align}
\label{eq_2norm_gram_err_general0}
&\min_{(S,G) \in {\cal R}} \|\K\|_2^2 \\
&= \min_{(S,G) \in {\cal R},\  \M \ \text{s.t.} \; \eqref{M_Lyap}} \begin{bmatrix} B \\ G \end{bmatrix}^T \begin{bmatrix} M_{11} & M_{12} \\ M_{12}^T & M_{22}\end{bmatrix} \begin{bmatrix} B \\ G \end{bmatrix}.\nonumber
\end{align}
Using the matrix notations  above, optimization problem \eqref{eq_2norm_gram_err_general0} can be written  in matrix form  explicitly as:
\begin{align}
\label{eq_2norm_gram_err_general_m0}
& \min_{(S,G,\M,\Pi)} \text{Trace}(\B_e^T \M \B_e) \\
& \text{s.t.:} \;\;  A \Pi + B L = \Pi S,  \; \sigma(S-GL) \subset \mathbb C^- \nonumber\\
& \qquad (S,G) \in \Gr, \;  \A_e^T \M + \M \A_e + \C_e^T \C_e =0. \nonumber
\end{align}
Note that,  the Sylvester equation $A \Pi + B L = \Pi S$ need \emph{not} be solved since, by \cite[Lemma 1] {i-astolfi-colaneri-SCL2014}, we can take $\Pi=VT$, with $V$ a certain Krylov projection  and $T$ some non-singular matrix. Therefore, we get the following simplified nonconvex optimization formulation for  Problem \ref{prob_optH2_general}:
\begin{align}
\label{eq_2norm_gram_err_general_m}
& \min_{(S,G,\M)} \text{Trace}(\B_e^T \M \B_e) \\
& \text{s.t.:} \;\; (S,G) \in {\cal R} \left( \Leftrightarrow \! \sigma(S \!-\! GL) \!\subset\!  \mathbb C^-\!,  (S \!-\! GL)_{ij} \!=\! 0 \; j \!\notin\! {\cal N}_i \right) \nonumber\\
&\qquad  \A_e^T \M + \M \A_e + \C_e^T \C_e =0,  \nonumber
\end{align}
In the rest of the paper we derive several numerical procedures for solving  the nonconvex problem \eqref{eq_2norm_gram_err_general_m}, whose optimal solution yields  a stable reduced order model of dimension $\nu$ of the linear network system, which preserves the topology of the network and minimizes the $H_2$-norm of the error system.


\section{Convex model reduction using block diagonal Gramians}
\label{sect_SDP}
\noindent The nonconvex problem \eqref{eq_2norm_gram_err_general_m} can be written equivalently in terms of matrix inequalities (semidefinite programming):
\begin{align}
\label{eq_2norm_gram_err_general_m22}
& \min_{(S, G) \in \Gr, \; \M \succeq 0 } \text{Trace}(\B_e^T \M \B_e) \\
& \text{s.t.:} \;\;   \A_e^T \M + \M \A_e + \C_e^T \C_e  \preceq 0,  \nonumber
\end{align}
where   $L$ is fixed a priori.  Clearly, semidefinte program (SDP) \eqref{eq_2norm_gram_err_general_m22} is not convex since it contains bilinear matrix inequalities (BMIs). However, our next result shows that  we can obtain a suboptimal solution through convex SDP using a simple assumption that the error  system admits a
block diagonal observability Gramian.  While diagonal Gramians have recently been exploited in the balanced truncation model reduction of positive systems \cite{BruRan:14}, the application of block diagonal Gramians on structured moment matching model reduction of general network systems is discussed in our paper.  

\begin{theorem}
\label{th:sdp}
If the  convex SDP relaxation:
\begin{align}
\label{eq_2norm_gram_err_general_m5}
& \min_{(X_{22},Y_{22},Z_{22},\Theta_{22}), M_{11} \succeq 0, M_{22} \succeq 0} \text{Trace} \left( B^T M_{11}B + X_{22} \right) \nonumber \\
& \text{s.t.}:  \Theta_{22}^T - L^T Z_{22}^T + \Theta_{22} - Z_{22}L +  (C\Pi)^T(C\Pi) \preceq Y_{22} \nonumber \\
&\qquad \begin{bmatrix}  X_{22} & Z_{22}^T \\ Z_{22} & M_{22}\end{bmatrix} \succeq 0, M_{22} \; \text{is block diagonal}  \\
& \qquad  \; ( \Theta_{22}  -  Z_{22} L)_{ij} =0 \;\; \forall i \in [N], j \notin {\cal N}_i \nonumber \\
&\qquad \begin{bmatrix} A^TM_{11} + M_{11}A +C^TC &  -C^T (C\Pi) \\ - (C\Pi)^T C & Y_{22} \end{bmatrix}  \preceq 0 \nonumber
\end{align}
has a solution, then we can recover a suboptimal solution of the model reduction Problem  \ref{prob_optH2_general}  expressed in terms of the SDP  problem \eqref{eq_2norm_gram_err_general_m22}  through the relations: 
\[ S= M_{22}^{-1} \Theta_{22}, \quad G= M_{22}^{-1} Z_{22},  \quad \M=\diag(M_{11}, M_{22}).   \] 
\end{theorem}

\begin{IEEEproof}
Using the block form of the Gramian $\M$, \eqref{eq_2norm_gram_err_general_m22} yields the equivalent SDP problem  \eqref{eq_2norm_gram_err_general_m3}.
\begin{figure*}[!htbp]
\begin{align}
\label{eq_2norm_gram_err_general_m3}
& \min_{(S,G)  \in \Gr, \; \M \succeq 0} \text{Trace} \left( B^T M_{11}B + B^TM_{12}G + G^TM_{12}^T B + G^TM_{22}G \right) \\
& \text{s.t.}:
\begin{bmatrix} A^TM_{11} + M_{11}A +C^TC & A^T M_{12} +M_{12} (S-GL) -C^T (C\Pi) \\
M_{12}^TA + (S-GL)^TM_{12}^T  - (C\Pi)^T C & (S-GL)^T M_{22} + M_{22}(S-GL) +  (C\Pi)^T (C\Pi) \end{bmatrix}  \preceq 0. \nonumber\\
\nonumber\\ \hline\nonumber
\end{align}
\end{figure*}
\noindent Note that problem \eqref{eq_2norm_gram_err_general_m3} is not convex since  if we assume $M_{12} \not =0$, then we cannot convexify the previous BMIs since we need to define $M_{12} G= Z_{12}$ and $M_{22} G= Z_{22}$ and require $\M \succeq 0$.  However, if we assume for the Gramian $\M$ the block $M_{12} =0$ and $M_{22}$  block diagonal, then problem \eqref{eq_2norm_gram_err_general_m3} can be recast as a convex SDP. More precisely, if we introduce additional variables, then we get the following SDP:
\begin{align}
\label{eq_2norm_gram_err_general_m4}
& \min_{(S,G,X_{22},Y_{22}), M_{11} \succeq 0, M_{22} \succeq 0} \text{Trace} \left( B^T M_{11}B + X_{22} \right) \\
& \text{s.t.}: X_{22} \succeq G^T M_{22} G,  (S-GL)_{ij} =0 \; \forall i \in [N], j \notin {\cal N}_i  \nonumber\\
&\qquad (S \!-GL)^T M_{22} + M_{22}(S \!-GL) +  \!(C\Pi)^T\!(C\Pi) \preceq Y_{22}, \nonumber \\
&\qquad \begin{bmatrix} A^TM_{11} + M_{11}A +C^TC &  -C^T (C\Pi) \\
- (C\Pi)^T C & Y_{22} \end{bmatrix}  \preceq 0.  \nonumber
\end{align}
Letting $Z_{22} = M_{22}G, \Theta_{22} = M_{22}S$ and using the Schur complement,  problem \eqref{eq_2norm_gram_err_general_m4} becomes the convex SDP \eqref{eq_2norm_gram_err_general_m5}. However, this change of variables  $M_{22}(S-GL) =  \Theta_{22} -  Z_{22} L$ is  in general not suitable when imposing the structured constraints $\Gr$ on $S-GL$.  Although the constraint on the parameterization $(S,G) \in \Gr$  is linear and thus convex, the corresponding constraint on $\Theta_{22}, Z_{22} $ and $M_{22}$ (i.e. $M_{22}^{-1} (\Theta_{22} -  Z_{22} L) \in \Gr$) is nonlinear and consequently nonconvex.  If we restrict the structure of $M_{22}$, assuming it is block diagonal with the block sizes compatible to those of the reduced  subsystems, i.e., $M_{22} = \text{diag}(M_{22}^1,\dots, M_{22}^N)$, with $M_{22}^i \in \mathbb R^{\nu_i \times \nu_i}$, the structured constraints are naturally guaranteed:
\begin{align}  
\label{eq:sdp_network}
M_{22}^{-1}(\Theta_{22} -  Z_{22}  L) \in \Gr  \;  \iff \;    \Theta_{22} -  Z_{22}  L \in \Gr. 
\end{align}
Note that the block diagonal assumption on $M_{22}$ is a sufficient condition for \eqref{eq:sdp_network}  given an arbitrary network structure $\Gr$. Moreover, we can recover a suboptimal solution of the original problem through the relations:  $G= M_{22}^{-1} Z_{22}, S= M_{22}^{-1} \Theta_{22}$ and  $ \M=\text{diag}(M_{11}, M_{22})$.  
\end{IEEEproof}

\noindent Clearly, \eqref{eq_2norm_gram_err_general_m5} is a suboptimal solution of the original SDP problem \eqref{eq_2norm_gram_err_general_m22}  since we restrict the Gramian matrix  $\M$ to have the   blocks $M_{12} =0$ and $M_{22}$ to be block diagonal. Hence, \eqref{eq_2norm_gram_err_general_m5} is  a convex SDP relaxation of the original problem \eqref{eq_2norm_gram_err_general_m22}.


\subsection{Sufficient conditions on block diagonal Gramians}
\noindent As we can see from the proof of Theorem \ref{th:sdp}, the block diagonal assumption on the observability Gramian is crucial. In this section we  derive sufficient conditions for the feasibility of the  SDP relaxation \eqref{eq_2norm_gram_err_general_m5}    in the case of general linear systems, i.e. sufficient conditions to guarantee that the error system admits a block diagonal observability Gramian.  For this we need the following result that holds for any two vectors $u$ and $v$:
\begin{align}
\label{eg:lmiineq}
u^T v + v^T u  \preceq u^T P^{-1} u + v^T P v \qquad \forall P \succ 0.
\end{align}

\noindent This inequality follows from the  relation $(u - P v)^T P^{-1} (u - P v) \succeq 0$ for $P \succ 0$.  We are interested in deriving sufficient conditions to guarantee that the SDP \eqref{eq_2norm_gram_err_general_m22}  admits a feasible triplet $(S,G,\M)$ with $\M$ of  block diagonal form and consequently the SDP relaxation \eqref{eq_2norm_gram_err_general_m5}  is well-defined.

\begin{theorem}
Given the stable minimal system \eqref{system} there exists a stable  reduced order model \eqref{redmod_CPi} such that the error system admits a block diagonal observability  Gramian $\M=\diag(M_{11}, M_{22})$, with $M_{22} = \text{diag}(M_{22}^1,\dots, M_{22}^N)$,  if the following conditions hold:
\begin{align}
\label{eq:sufcond}
& A^T M_{11}  \!+\! M_{11}A  \!+\! C^TC  \!+\!  C^T (C\Pi) P^{-1}\!(C\Pi)^T C \preceq 0   \\
&(S-GL)^T M_{22} + M_{22} (S-GL) + (C\Pi)^T (C\Pi)  +  P  \preceq 0  \nonumber
\end{align}
for some matrices   $M_{11}, M_{22}, P \succeq 0$ and $(S,G) \in {\cal G}$.
\end{theorem}

\begin{IEEEproof}
Note that the feasible set of SDP  \eqref{eq_2norm_gram_err_general_m22}  is nonempty if   $\M \succeq 0$, $(S,G) \in {\cal G}$ and the following inequality holds:
\begin{align}
\label{eq:bd}
& \begin{bmatrix}  x \\ \xi \end{bmatrix}^T \! \left( \! \begin{bmatrix} A & 0 \\ 0 & S \!-\!GL \end{bmatrix}^T \! \M  \!+\! \M \! \begin{bmatrix} A & 0 \\ 0 & S \!-\! GL \end{bmatrix} \right) \!\! \begin{bmatrix}  x \\ \xi \end{bmatrix}  \\ 
& \qquad  + \begin{bmatrix}  x \\ \xi \end{bmatrix}^T \! \left(  \begin{bmatrix} C^TC & -C^T C\Pi \\ -(C \Pi)^T C & (C\Pi)^T C\Pi  \end{bmatrix}  \right) \! \begin{bmatrix}  x \\ \xi \end{bmatrix} \leq 0 \quad \forall x,\xi,   \nonumber
\end{align}
which, using that $\M=\diag(M_{11}, M_{22})$,  is equivalent to
\begin{align*}
& x^T(A^T M_{11} + M_{11}A + C^TC)x \\
& \quad + \xi^T((S-GL)^T M_{22} + M_{22} (S-GL) + (C\Pi)^T (C\Pi))\xi \\
&\quad - x^T C^T (C\Pi) \xi - \xi^T (C\Pi)^T C x \leq 0 \quad \forall x,\xi.
\end{align*}
Using now  \eqref{eg:lmiineq} in the last term we get that:
\[   - x^T\! C^T \! (C\Pi) \xi - \xi^T \! (C\Pi)^T \! C x  \! \leq \!  x^T \! C^T \! (C\Pi) P^{\!-1} \! (C\Pi)^T \! C x + \xi^T \! P  \xi   \]
for all $x,\xi$.   Consequently, if    the  inequality:
\begin{align*}
& x^T(A^T M_{11} + M_{11}A + C^TC +  C^T  (C\Pi)  P^{-1} (C\Pi)^T C)x \\
& \!+\! \xi^T \! ((S \!-\! GL)^T \! M_{22} \!+\! M_{22} (S \!-\!GL) \!+\! (C\Pi)^T\!(C\Pi) \!+\! P) \xi \!\leq\! 0 
\end{align*}
holds for all $x,\xi$,  then  \eqref{eq:bd}  also holds. This proves  the sufficient conditions \eqref{eq:sufcond}.
\end{IEEEproof}

\noindent This theorem provides  sufficient conditions  and  a  procedure for constructing a reduced order network model for which  the corresponding  error system admits a bock diagonal observability Gramian. Indeed, let us, for example, fix $L, M_{22}= I_\nu$, and some $ P \succ 0$. Then, the existence of a solution $(S,G,\Pi)$ satisfying $(S, G) \in {\cal G}$ of the system of equations:  
\[  (S-GL)^T \!+\!  (S-GL)  \!=\! -\!\left(\!(C\Pi)^T\!(C\Pi) \!+\! \!P \right)\!,    A \Pi + BL \!=\! \Pi S \]
together with the existence of an $M_{11} \succ 0$ satisfying  $A^T M_{11} + M_{11}A + C^TC +  C^T (C\Pi) P^{-1} (C\Pi)^T C \preceq 0$ guarantee that we have a reduced order model for which the corresponding observability Gramian of the error system is   block diagonal (i.e. $M_{12} =0$ and $M_{22} = I_\nu$ is also block diagonal). For example, we can fix  matrices $S$ , $L$ and $G$ such that $S-GL$ is stable and has the same network structure as $A$ (just take  stable diagonal matrix $S$  and $G=0$). Let $\Pi$ be the solution of  $A \Pi + BL = \Pi S$ and define $P = -(S-GL)^T - (S-GL) - (C\Pi)^T (C\Pi)$. If the resulting  $P \succ 0$ and if there exists  $M_{11} \succ 0$ satisfying  $A^T M_{11} + M_{11}A + C^TC +  C^T (C\Pi) P^{-1} (C\Pi)^T C \preceq 0$, then we obtain a stable reduced  model preserving the network structure and for which the corresponding error system admits a block diagonal  Gramian $\M = \text{diag}(M_{22}, I_\nu)$.  

\noindent From our best knowledge, the most common dynamical systems that admit block diagonal Gramians are the positive systems. The system matrices for these systems satisfy:  all off-diagonal elements of the matrix $A$ and all the entries of the matrices $B$ and $C$ are non-negative.  Positive systems occur in modelling of applications with special structures from, e.g.,  biomedicine, economics,  networks (see \cite{SooAnd:14,ReiVir:09}).


\section{ Nonconvex model reduction using projected gradient}
\label{sect_gradient}
\noindent  The  SDP approximation \eqref{eq_2norm_gram_err_general_m5} offers a convex way to solve the nonconvex problem  \eqref{eq_2norm_gram_err_general_m} of Problem 1 at the cost of some performance loss in general. We can also derive a numerical approach based on projected   gradient method to solve \eqref{eq_2norm_gram_err_general_m}. Here, our idea is to use a partial minimization approach (see Appendix \ref{note_on_parmin}) to  \eqref{eq_2norm_gram_err_general_m}  leading to a smooth reformulation, and apply the gradient projection method to get a (locally) optimal solution for \eqref{eq_2norm_gram_err_general_m}. More precisely,  consider the nonconvex optimization problem \eqref{eq_2norm_gram_err_general_m}, where  $L$  fixed a priori and  $\A_e = \A_e(S,G)$.  Then, the following partial minimization holds for  problem \eqref{eq_2norm_gram_err_general_m}: 
\begin{align*}
&\eqref{eq_2norm_gram_err_general_m} = \min_{(S,G) \in {\cal R}}  \left(  \min_{\M: \A_e^T \M + \M \A_e + \C_e^T \C_e =0} \text{Trace}(\B_e^T \M \B_e)) \right).
\end{align*}
However, if $S-GL$ and  $A$ are stable,  then there  exists unique  $\M =\M(S,G)$ positive semidefinte   solution of the Lyapunov equation:
\[  \A_e^T \M + \M \A _e+ \C_e^T \C_e =0.\] 
Hence, for any pair $(S,G)$ stable, the partial minimization in $\M$ leads to an optimal value $f(S,G) = \min_{\M: \A_e^T \M + \M \A_e + \C_e^T \C_e =0} \text{Trace}(\B_e^T \M \B_e))$, which can be written explicitly as:
\[  f(S,G) = \text{Trace} \left( \begin{bmatrix} B \\ G \end{bmatrix}^T  \M(S,G) \begin{bmatrix} B \\ G \end{bmatrix} \right),  \] 
where $\M(S,G)$ is the unique positive semidefinite  solution of the Lyapunov equation:
\begin{align}
\label{lyap_eq_11}
&\begin{bmatrix} A & 0 \\ 0 & S \!-\! GL \end{bmatrix}^T \!\!\!\! \M   \!+\!  \M \!\!  \begin{bmatrix} A & 0 \\ 0 & S \!-\! GL \end{bmatrix} \!+\! \begin{bmatrix} C^TC & -C^T \! C\Pi \\ 
-(C\Pi)^T \! C & (C\Pi)^T \! (C\Pi)  \end{bmatrix} \nonumber \\ 
&  = 0.
\end{align}
Thus,  we get the following equivalent reformulation for \eqref{eq_2norm_gram_err_general_m}:
\begin{align}
\label{eq_2norm_gram_err_general_m2}
& \min_{(S,G)} f(S,G) \;\;  \left (:= \text{Trace} \left( \begin{bmatrix} B \\ G \end{bmatrix}^T \!\! \M(S,G) \begin{bmatrix} B \\ G \end{bmatrix} \right) \right) \\
& \text{s.t.}: \;\; (S,G) \in {\cal R}  \quad \text{and} \quad \eqref{lyap_eq_11}. \nonumber
\end{align}
For solving the equivalent   problem \eqref{eq_2norm_gram_err_general_m2} we can  apply any first- or second-order optimization method. Hence, we need to compute the gradient and eventually  the Hessian of the objective function $f$. In the sequel, we show that we can compute the gradient of the objective function of  \eqref{eq_2norm_gram_err_general_m2} solving two Lyapunov equations.  Indeed, since $\text{Trace}(MN) = \text{Trace}(NM)$ for any matrices $M, N$ of compatible sizes,  the  objective function of \eqref{eq_2norm_gram_err_general_m2} becomes:
\begin{align*}  
f(S,G) &  = \text{Trace} \left(  \M(S,G)   {\cal B}(S,G)  \right), \;  {\cal B}(S,G)  =  \begin{bmatrix} B \\ G \end{bmatrix} \begin{bmatrix} B \\ G \end{bmatrix}^T.
\end{align*}

\begin{theorem}
\label{th:pm}
The  objective function  $f$ of optimization problem \eqref{eq_2norm_gram_err_general_m2}  is differentiable on the set of stable matrices ${\cal D} = \left\{ (S, G): \;  \sigma(S-GL) \subset \mathbb C^-  \right\}$  and the gradient of $f$ at any pair of matrices  $(S,G) \in {\cal D}$ is given by $\nabla f(S,G) = [\nabla_S f(S,G)  \; \nabla_G f(S,G)] \in  {\mathbb R}^{\nu \times (\nu +m)}$ with: 
\begin{align}
\label{eq_gradient1}
& \nabla_S f(S,\!G)   \!=  2   \left[ M_{12}^T(S,\!G)  W_{12}(S,\!G)  \!+\!  M_{22} (S,\!G) W_{22}(S,\!G) \right] \nonumber \\
& \nabla_G f(S,\!G)  \!= 2   \left[M_{12}^T(S,G) B +  M_{22}(S,G) \ G \right. \\
& \qquad \left.   - M_{12}^T(S,G)  W_{12}(S,\!G) L^T  \!-\! M_{22} (S,G) W_{22}(S,G) L^T   \right],  \nonumber
\end{align}
where $\M(S,G) $ solves the Lyapunov equation \eqref{lyap_eq_11} and   $\W(S,G) $ solves the Lyapunov equation from below:
\begin{equation}
\label{WX_lyap}
\begin{bmatrix} A & 0 \\ 0 & S \!-\! GL \end{bmatrix} \!  \W(S,G) + \W(S,G) \!  \begin{bmatrix} A & 0 \\ 0 & S \!-\! GL\end{bmatrix}^T \!+   {\cal B}(S,G) \!=\! 0.
\end{equation}
\end{theorem}

\begin{IEEEproof}
To compute the gradient $\nabla f(S,G)$, we write the derivative  $f'(S,G) \opd{(S,G)}$ for some $\opd{(S,G)} \in  \mathbb R^{\nu \times (\nu+m) }$ in gradient form using the trace. We introduce the gradient as:
\begin{align*} 
f'(S,G)  & \opd{(S,G)}  = \text{Trace} \left( \nabla f(S,G)^T \opd{(S,G)}  \right) \\
& = \text{Trace} \left( \nabla_S f(S,G)^T \opd{S} +  \nabla_G f(S,G)^T \opd{G} \right). 
\end{align*}
Then, we have:
\begin{align*}
&f'(S,G)  \opd{(S,G)}\\  
&\quad =  \text{Trace} \left(     \M'(S,G) {\cal B}(S,G)   + \M(S,G) {\cal B}'(S,G)   \right).
\end{align*}
We  compute separately the two terms in the above expression.  Let us define:  $$\Phi(S,G,\M)=\begin{bmatrix} A & 0 \\ 0 & S-GL \end{bmatrix}^T \M  + \M  \begin{bmatrix} A & 0 \\ 0 & S-GL \end{bmatrix}.$$ 
Since $(S,G) \in {\cal D}$ and ${\cal D}$ is an open set, then  $ \Phi_{\M}(S,G,\M) \opd \M$ given by:
\begin{align*}
 \Phi_{\M}(S,G,\M) \opd \M &=\begin{bmatrix} A & 0 \\ 0 & S \!-\! GL \end{bmatrix}^T \!\! \opd\M  + \opd\M  \begin{bmatrix} A & 0 \\ 0 & S \!-\! GL \end{bmatrix}
\end{align*}
is surjective and also we have:
\begin{align*} 
& \Phi_{(S,G)} (S,G,\M)  \opd{(S,G)}\\ 
&= \begin{bmatrix} 0 & 0 \\ 0 & \opd{S}  - \opd{G} \  L \end{bmatrix}^T \M + \M \begin{bmatrix} 0 & 0 \\\ 0 & \opd{S}  - \opd{G} \ L  \end{bmatrix}.
\end{align*}
Since $\Phi(S,G,\M) +{\cal C}_e^T {\cal C}_e =0$,   the Implicit Function Theorem yields the differentiability of $\M(S,G)$ and the relation:
\begin{align}
\label{eq_TFI1}
& \begin{bmatrix} A & 0 \\ 0 & S-GL  \end{bmatrix}^T \M'(S,G)  + \M'(S,G) \begin{bmatrix} A & 0 \\ 0 & S-GL \end{bmatrix}  \\
&+ \! \begin{bmatrix} 0 & 0 \\ 0 & \opd{S}  \!-\! \opd{G}   L  \end{bmatrix}^T \!\!\! \M(S,G) + \M(S,G)\! \begin{bmatrix} 0 & 0 \\  0 & \opd{S}  \!-\! \opd{G}   L \end{bmatrix} \!=\! 0.  \nonumber
\end{align}

\noindent Moreover, by \eqref{W_Lyap} the  Gramian $\W(S,G)$ is the unique solution of the Lyapunov equation \eqref{WX_lyap}. Subtracting \eqref{eq_TFI1} multiplied by $\W(S,G)$ to the left from \eqref{WX_lyap} multiplied by $\M'(S,G) $ to the right,  taking the trace,  and reducing the appropriate terms, we get the relation:
\begin{align}  
\label{eq:df1}
&\text{Trace} \left( \M'(S,G)   {\cal B}(S,G)  \right) \\  
&=  \text{Trace}  \left(  \W(S,G) \begin{bmatrix} 0 & 0 \\ 0 &  \opd{S}  - \opd{G}   L  \end{bmatrix}^T  \M(S,G)  \right. \nonumber \\ 
& \qquad \qquad \qquad \left. + \M(S,G) \begin{bmatrix} 0 & 0 \\ 0 &  \opd{S}  - \opd{G}   L  \end{bmatrix}  \W(S,G) \right). \nonumber
\end{align}
Similarly, for the second term  we get: 
\begin{align}  
\label{eq:df2}
&\text{Trace} \left(  M(S,G)  {\cal B}'(S,G)   \right) \\
& =  \text{Trace} \left( \M(S,G)  \begin{bmatrix} 0 & B   \opd{G}^T   \\    \opd{G} B^T &  \opd{G} \ G^T + G   \opd{G}^T \end{bmatrix}   \right).  \nonumber
\end{align}
Hence, combining  \eqref{eq:df1} and \eqref{eq:df2},  using the block structure of $\W$ and $\M$, and the definition of trace, we obtain easily the closed form expression for the gradient from \eqref{eq_gradient1}.
\end{IEEEproof}


\noindent The  previous theorem  also yields the necessary optimality condition for the model reduction  optimization problem \eqref{eq_2norm_gram_err_general_m2}:
\begin{lemma}
\label{lem_optim}
Let  the block presentations $S=(S_{ij})_{i,j=1:N}$,  $G=(G_i)_{i=1:N}$  and $L=(L_i)_{i=1:N}$, with $L$ fixed.  If   $(S^*,G^*) \in {\cal D}$  solves the optimization problem   \eqref{eq_2norm_gram_err_general_m2}, then 
\[   \nabla_{S_{ij}}  f(S^*,G^*)  - \nabla_{G_i}  f(S^*,G^*) \cdot  L_j = 0 \quad \forall i \in [N],  j \in {\mathcal N}_i,  \] 
where  the expressions of $\nabla_S f$  and   $\nabla_G f $  are given in \eqref{eq_gradient1}. 
\end{lemma}

\noindent We can replace the open set ${\cal D}$ with any sublevel set:
\[  {\cal L}(S_0,G_0) = \{ (S,G) \in  {\cal D}: \; f(S,G) \leq f(S_0,G_0)  \},  \]
where $(S_0,G_0) \in  {\cal D}$ is any pair of  initial stable reduced order system matrices. By arguments as in \cite{Toi:85} we can show that ${\cal L}(S_0,G_0)$ is a compact set.   Then,  the theorem of Weierstrass implies that for any  given matrix $(S_0,G_0) \in  {\cal D}$,  the model reduction Problem  \ref{prob_optH2_general}  given by  optimization formulation \eqref{eq_2norm_gram_err_general_m2} has a global minimum  in the sublevel set ${\cal L}(S_0,G_0)$. We can also show that the gradient $\nabla f(S,G)$ is Lipschitz continuous on the compact sublevel  set   ${\cal L}(S_0,G_0)$.  We briefly sketch the proof of this statement. First we observe that $\M(S,G)$ and $\W(S,G)$ are continuous functions, since they are solutions of some algebraic linear systems. Moreover, there exists finite $\ell_M > 0$ such that for all $(S,G), (S',G') \in   {\cal L}(S_0,G_0)$:
\[   \|  \M(S,G) -   \M(S',G')  \| \leq \ell_M \|  (S,G)  - (S',G') \|.  \]
Then, using the expression of  $\nabla f(S,G)$, compactness of  ${\cal L}(S_0,G_0)$, continuity of $\M(S,G)$ and $\W(S,G)$, and the previous relation we conclude that there exists $\ell_f >0$ such that for all $(S,G), (S',G') \in   {\cal L}(S_0,G_0)$:
\[  \|  \nabla f(S,G) -   \nabla f (S',G')  \| \leq \ell_f  \|  (S,G)  - (S',G')\|.  \]
This property of the gradient is useful when analyzing the convergence behavior of the projected gradient algorithm  we propose below for solving optimization problem  \eqref{eq_2norm_gram_err_general_m2}.


\subsection{Projected gradient method}\label{sect_proj_grad}
\noindent We have proved that the nonconvex optimization problem \eqref{eq_2norm_gram_err_general_m2} has differentiable objective function  and its gradient is given by \eqref{eq_gradient1}. Moreover, the gradient is smooth  (i.e. Lipschitz continuous) on any compact set. Then, we can apply the projected gradient method for obtaining a (local) optimal solution of  \eqref{eq_2norm_gram_err_general_m2}. Starting from  an initial stable matrix pair  satisfying the structured constraints,  $(S_0, G_0)  \in {\cal R}$,  we consider the following update rule:
\begin{align}  
\label{eq:pgrad}
(S_{k+1},G_{k+1})  =  (S_k,G_k) - \alpha_k  \Pi_{\cal G} \left( \nabla f(S_k,G_k) \right), 
\end{align}
where the stepsize $\alpha_k$ can be chosen by a backtracking procedure or constant in the interval  $(0, \ 2/\ell_f)$ (where $\ell_f$ denotes the Lpschitz constant of the gradient) \cite{Nes:04}. Here  $ \Pi_{\cal G} \left( \nabla f \right)$ denotes the projection of the gradient of the objective function $\nabla f$ onto the convex set ${\cal G}$ describing the network structure.   Note that the projection of  $\nabla f$ onto the convex set ${\cal G}$ is straightforward and computationally cheap: using the expressions of the gradient from \eqref{eq_gradient1}, we only need to set the blocks of the gradient for all  $i \in [N]$ and $j \notin  {\cal N}_i$ as 
\[ (\nabla_S f)_{ij} = (\nabla_G f \cdot  L)_{ij} \quad \iff \quad \nabla_{S_{ij}} f = \nabla_{G_i} f \cdot  L_j. \]  

\noindent Based on this update law and with these choices for the stepsize, using the Lipschitz gradient property for the objective function we can easily prove that the sequence of value functions $f(S_k,G_k)$ is nonincreasing \cite{Nes:04}:
\[  f(S_{k+1},G_{k+1}) \leq f(S_k,G_k) - \Delta \cdot \| (S_{k+1},G_{k+1}) - (S_k,G_k)  \|^2  \]
for some constant $\Delta>0$ for all $k \geq 0$.  Therefore all the iterates remain in the compact sublevel set ${\cal L}(S_0,G_0)$.   Moreover,  we define gradient mapping of $f$  at $(S_k,G_k)$  \cite{Nes:04}: 
\[ \Gamma_{\cal G} (S_k,\!G_k) \!=\! \alpha_k^{-1} (\!(S_k,\!G_k) - (S_{k \!+\! 1},G_{k \!+\!1})\!) \!=\!   \Pi_{\cal G} \! \left( \nabla f(S_k,\!G_k) \right)\!,  \]    
which, according to Lemma \ref{lem_optim}, represents the natural measure of optimality for the constrained problem  \eqref{eq_2norm_gram_err_general_m2}.  Since $f$ is bounded from below by zero, then for any positive integer  $k>0$ it is straightforward to prove from the previous descent inequality   the following global convergence rate for the gradient mapping: 
\[   \min_{t=0:k}  \| \Gamma_{\cal G} (S_t,G_t)  \|^2   \leq    \frac{f(S_0,G_0) - f^*}{\Delta \cdot  (\min_{t=0:k} \alpha_t)  \cdot k}  \quad \forall k > 0, \]
where $f^*$ is the optimal value of problem \eqref{eq_2norm_gram_err_general_m2}. Moreover,   under some mild assumptions, such as the Hessian of $f$ at a local minimum is positive definite and  bounded, then starting sufficiently close to this local optimum the gradient iteration converges linearly to this solution \cite{Nes:04}. Therefore, the speed of convergence of this iterative process depends on the starting point.  Hence, to obtain a good  initial stabilizing pair of matrices $(S_0,G_0)$ satisfying the network conditions $\Gr$, we can  solve the structured SDP problem \eqref{eq_2norm_gram_err_general_m5} and initialize with its solution.    Note that for computing the gradient we first need  to find the Gramians  $\M_k$ and $\W_k$,   solutions of the Lyapunov equations  \eqref{lyap_eq_11} and  \eqref{WX_lyap} in $(S_k, G_k)$. The solvability of $\M_k$ and $\W_k$, unique positive semidefinite  solutions of \eqref{lyap_eq_11} and  \eqref{WX_lyap},  is implied by the stability of the error matrix $\A_e(S_k,G_k)$. 

\noindent Hence, we consider the following algorithmic procedure:
\begin{alg}[$H_2$ optimal network reduction algorithm]
\begin{enumerate}[wide,nosep]
\item Let $(S_0,G_0) \in \mathcal{G}$ (e.g., $(S_0,G_0)$ - solution of SDP \eqref{eq_2norm_gram_err_general_m5}).
\item Perform update \eqref{eq:pgrad} until $\| \Gamma_{\cal G} (S_k,G_k)  \|\leq\epsilon$.
\end{enumerate}
\end{alg}
Note that the Gramians  $\M_k$ and $\W_k$ are in general dense, which means the block diagonal assumption from Section \ref{sect_SDP}  is relaxed during each iteration of our projected gradient method. In fact, the block diagonal Gramian can be viewed as an intermediate step between a diagonal Gramian in \cite{BruRan:14}  on positive linear systems and a full one as provided by our gradient method on general linear network systems. 


\section{Extensions}
\label{sect_ext}
\noindent Note that our approach is  general and  flexible, allowing to tackle network systems with even more structure. For example, we can consider that coupling among subsystems is   through both,  states and  inputs, i.e. we  modify the dynamics \eqref{mod3}  as:
\begin{align}
\label{mod3B}
\dot x_i & =  \sum_{j \in {\mathcal N}_i^x}  A_{ij} x_{j} +  \sum_{j \in {\mathcal N}_i^u} B_{ij} u_{j}   \quad  \forall i \in [N],
\end{align}
where  $x_i$ and $u_i$ represent the states and inputs of the $i$th subsystem. Our optimization-based model reduction framework allows to easily  incorporate the additional structured constraints, inherited from $B$ (i.e. $B_{ij}=0 \; \forall i \in [N], j \notin {\cal N}_i^u$),  on the matrix $G$.  For example in the convex SDP  \eqref{eq_2norm_gram_err_general_m5} we just need to add the additional convex constraint $ (Z_{22})_{ij} = 0  \; \forall i \in [N], j \notin {\cal N}_i^u$. Similarly, in the projected gradient method we just need to set the corresponding block components of the gradient $\nabla_G f$ to zero,  i.e. $\nabla_{G_{ij}} f = 0  \; \forall i \in [N], j \notin {\cal N}_i^u$. Similarly, for positive networks we can easily incorporate positivity constraints on the system matrices of the reduced model.


\section{Illustrative examples}
\label{sect_exmp}
\noindent In this section, we illustrate the efficiency of the proposed methods numerically. In particular, we compute and compare  reduced order models for two network examples  achieving (possibly) the minimum  $H_2$ norm approximation.


\subsection{Random positive network system}
\label{sect_ex_net_rand}
\begin{figure}[!htbp]
\centering
\includegraphics[width=.5\textwidth]{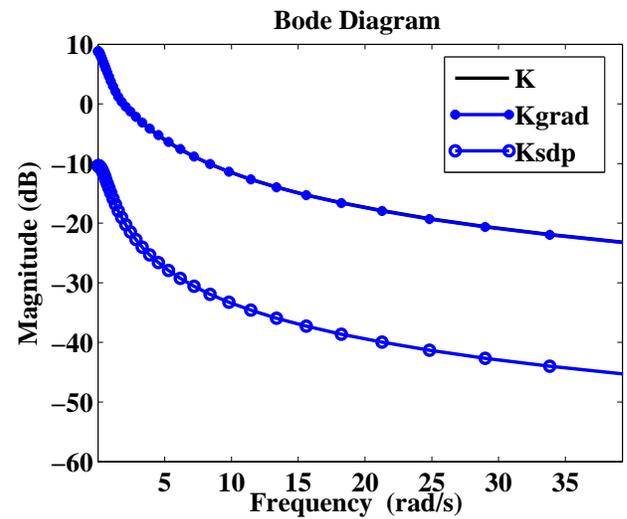}
\caption{Bode magnitude plots of the given $12$th order network (black) vrs. the $4$th order $H_2$ SDP and  optimal gradient newtork aproximations (blue).}
\label{resultat_randnet}
\end{figure}

\noindent We consider a stable $12$th order linear positive  network system as in \eqref{mod3} with matrices $A, B$ and $C$ described in Appendix \ref{sect_ap_par}, generated randomly in the interval $(-5, 1)$ such that $A$ is stable and positive and satisfying the interconnection map from Figure  \ref{fig:net_systems}. Here, we select $n_i = 3$, $m=1$ and $p=1$ for all $N=4$ subsystems. Hence, $A \in \mathbb{R}^{12 \times 12}$.  We compute  an $H_2$ optimal reduced order  network, with  $\nu_i=1$ for all $i=1:4$, with the same network structure. For the initialization of the algorithm we consider the solution   $(S_0,G_0)$  of the structured SDP problem \eqref{eq_2norm_gram_err_general_m5}, with $ L=[0 \;  0 \;  0 \; 1]$ fixed a priori.   The output responses of $H_2$ SDP $\widehat K_{\text{sdp}}(s)$ and  $H_2$ optimal obtained with projected gradient $\widehat K_{\text{grad}}(s)$ are displayed in  Figure \ref{resultat_randnet}. We observe that the $H_2$ optimal output response, corresponding to  projected gradient $\widehat K_{\text{grad}}(s)$, is almost identical with the response of $K(s)$ in the frequency range we considered  $10^{-6}: 10^2$. 

\noindent In Table \ref{tab_comp}, we compare the $H_2$-norm of the errors yielded by the proposed SDP and projected gradient methods versus the positive preserving projections method in \cite{SooAnd:14} and the positivity preserving balanced truncation (BT) in \cite{ReiVir:09}. Note that our projected  gradient method performs best, that is  the  $H_2$ error  of the $4$th order network decreased from $2.813$ (for SDP) to $5.075\cdot 10^{-3}$ (for projected gradient).  Moreover, the constraints (iv) in Problem \ref{prob_optH2_general} are satisfied by the optimal model.

\begin{table}[!htbp]
\centering
\begin{tabular}{c|c|c|c|c}
Method &  SDP & Gradient & Projection \cite{SooAnd:14} & Positive BT \cite{ReiVir:09} \\
\hline
$H_2$ error & 2.813 & $5.075\cdot 10^{-3}$ & 1.9103 & 1.7689
\end{tabular}
\caption{$H_2$ error norm of four reduction methods. }
\label{tab_comp}
\end{table}


\subsection{Power network system}
\noindent Consider a power system split into $N$ control areas consisting of a generator and a load, with tie-lines providing interconnections between them, as described in \cite{Ven:06}. For the area $i$ the simplified model is given by the differential equations:
\begin{align}\label{eq_area}
\frac{\opd\Delta\omega_i}{\opd t} &= -\frac{D_i}{M_i^a}\Delta\omega_i-\frac{1}{M_i^a}\Delta P_{m_i}-\frac{1}{M_i^a}\Delta P_{\rm tie}^{ij}, \nonumber\\
\frac{\opd\Delta P_{m_i}}{\opd t} &= -\frac{1}{T_{{\rm CH}_i}}\Delta P_{m_i}+\frac{1}{T_{{\rm CH}_i}}\Delta P_{v_i}, \\
\frac{\opd\Delta P_{v_i}}{\opd t} &= -\frac{1}{R_i^fT_{G_i}}\Delta\omega_i-\frac{1}{T_{G_i}}\Delta P_{v_i}+\frac{1}{T_{G_i}}\Delta P_{{\rm ref}_i}. \nonumber
\end{align}
The interconnection to the control area $j\ne i$, is made through the tie-line with the power flow modeled by the equation:
\begin{equation}\label{eq_tie}
\frac{\opd\Delta P_{\rm tie}^{ij}}{\opd t} = T_{ij}(\Delta\omega_i-\Delta\omega_j), 
\end{equation}
with $\Delta P_{\rm tie}^{ij}=-\Delta P_{\rm tie}^{ji}$. The notation $\Delta$ indicates the deviation from the steady-state of the variable, e.g., $\Delta\omega$ is the deviation of the angular speed from the nominal operating value. The variable are defined as follows:
\begin{itemize}[wide]
\item $\omega$ is the angular speed of the rotor,
\item $M^a$ is the angular momentum,
\item $D$ is the ratio between the percentage change in load and the percentage change in frequency,
\item $P_m$ is the mechanical power,
\item $T_{\rm CH}$ is the charging time constant,
\item $P_v$ is the steam valve position,
\item $P_{\rm ref}$ is the load reference setpoint,
\item $R^f$ is the ratio between the percentage change in frequency and the percentage change in unit output,
\item $T_G$ is the governor time constant,
\item $P_{\rm tie}^{ij}$ is the tie-line power flow between areas $i$ and $j$,
\item $T_{ij}$ is stiffness coefficient  of the tie-line (between area $i$ and area $j$).
\end{itemize}

\begin{figure}[!htbp]
\centering
\includegraphics[width=.5\textwidth]{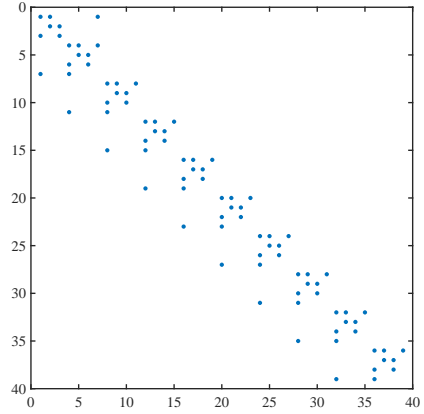}
\caption{Sparsity pattern of the $A$ matrix of a power network with $N=10$ areas described by equations \eqref{eq_area} and \eqref{eq_tie}}
\label{fig_powNetSpy}
\end{figure}

\noindent Considering a chain of $N$ areas, equations \eqref{eq_area} and \eqref{eq_tie} yield a $4N-1$th order model of the power network of the form \eqref{mod3B} with the states: 
\[ x_1=\begin{bmatrix}\Delta\omega_1 \\ \Delta P_{m_1} \\ \Delta P_{v_1}\end{bmatrix},\quad x_i=\begin{bmatrix} \Delta\omega_i \\ \Delta P_{m_i} \\ \Delta P_{v_i} \\ \Delta P_{\rm tie}^{ij} \end{bmatrix} \quad  i=2:N. \] 
The input of each area is $u_i = \Delta P_{{\rm ref}_i}$, while as measured outputs we consider the angular speed deviations $\Delta\omega_i$ of each area  $i=1:N$.  Note that the matrix $A$ of the power network model has almost block bidiagonal form as seen in Figure \ref{fig_powNetSpy}.

\begin{figure}[!htbp]
\centering
\includegraphics[width=.45\textwidth]{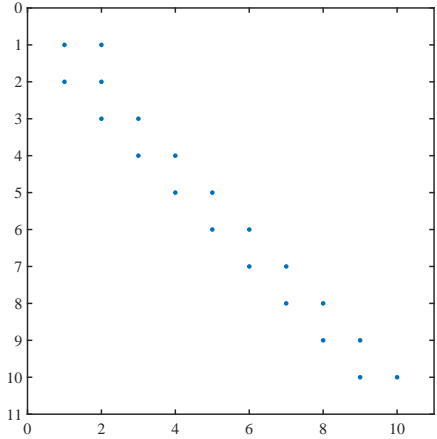}
\caption{Sparsity pattern of  $F=S-GL$ matrix of the $N$th order approximation of a power network with  $N=10$ areas. }
\label{fig_approxSpy}
\end{figure}

\noindent  The values of the parameters are chosen for each of the $N$ areas, randomly in the interval between the lowest and the highest possible physical values extracted from \cite{Ven:06}:
\begin{itemize}
\item $D_i\in [0.255, 80]$,
\item $M_a^i, T_{{\rm CH}_i}\in [1,5]$,
\item $R^f_i\in [0.03,0.07]$,
\item $T_{G_i}\in [4,10]$,
\item $T_{ij}\in [1.5,2.5]$,
\end{itemize}
with $i,j=1:N$. With these values the resulting model of the power network is \emph{stable}. For $N=4$, the matrices of the linear network model are given by \eqref{eq_ABC_ex_net_rand} and by \eqref{eq:BC}.
We aim at selecting $N$ interpolation points, that is the order of each reduced subsystem is  $\nu_i=1$ for all $i=1:N$, and accordingly the matrices $S, L, G \in\mathbb R^{N\times N}$.  For solving the model reduction Problem \ref{prob_optH2_general} for the power network model from above with a chain of $N$ areas we use the projected gradient algorithm from Section \ref{sect_proj_grad} initialized with the SDP solution \eqref{eq_2norm_gram_err_general_m5}, for $L=I_N$ fixed a priori.  This approach yields a $N$th order $H_2$ optimal approximation of the power network of the form \eqref{red_mod_F_block} redering an optimal $H_2$ norm of the approximation error. Note that the projected gradient method preserves the (block) bidiagonal network topology of $S-GL$ as Figure \ref{fig_approxSpy} shows, as well as the \emph{stability} of the network. 

\begin{figure}[!htbp]
\centering
\includegraphics[width=.5\textwidth]{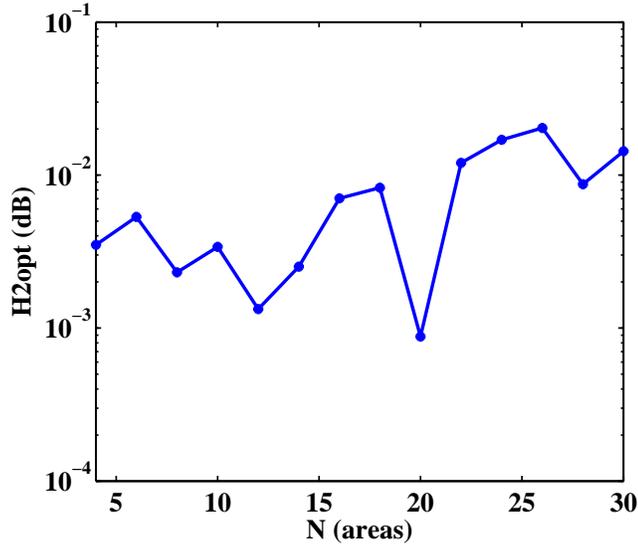}
\caption{Plot of $H_2$ optimal norm of approximation error for $N=4:30$.}
\label{fig_H2vN}
\end{figure}
\noindent We plot in Figure \ref{fig_H2vN}  the optimal $H_2$ norm of the error approximation yielded by the projected gradient method  for a number of $N$ areas ranging from  $4:30$.  The variation of the optimal $H_2$ norm is consistent with the increase in $N$, that is the $N$th order approximation of the $4N-1$ model is not necessarily more accurate as $N$ increases.


\section{Conclusions}
\noindent In this paper we have studied the  model order  reduction for linear network systems. Using moment matching techniques, we have developed an optimization framework to compute  parameterized reduced order stable  models  achieving moment matching, minimizing the $H_2$ norm of the error system, and preserving the structure of the to-be-reduced model of the network.  For this, we have  proposed two numerical procedures, based on SDP and projected gradient,  for finding the (optimal) reduced order model of the network.  Preliminary numerical simulations have confirmed the efficiency of our approach.

\begin{figure*}[!htbp]
\begin{align}\label{eq_ABC_ex_net_rand}
A\! &=\! \left[\!\begin{array}{cccccccccccc}
 -4.6000&1.0000&0&0.1000&0.0500&0.0200&0.0100&0&0&0&0&0 \\
  0&-4.6000&1.0000&0.1000&0.2100&0.1000&0&0.0500&0&0&0&0 \\
  0.4714&0.1953&-4.2667&0&0.0900&0.0800&0&0&0.0200&0&0&0 \\
  0.1000&0&0.3300&-4.2300&0&0.2000&0.1000&0&0&0&0&0 \\
  0.2000&0.0200&0&0.1000&-4.5980&0&0&0.0200&0&0&0&0 \\
  0.1000&0&0.1600&0&0.1210&-7.1900&0&0&0.0500&0&0&0 \\
  0.5324&0.3365&0.4039&0.5527&0.2431&0.9357&-3.3242&0.0012&0.6253&0.2874&0.7624&0.6455 \\
  0.7165&0.1877&0.5486&0.2748&0.1542&0.8187&0.3604&-3.1836&0.5431&0.5017&0.5761&0.1232 \\
  0.1793&0.3219&0.0487&0.2415&0.9564&0.7283&0.1888&0.6996&-3.0610&0.7615&0.7477&0.5044 \\
  0&0&0&0&0&0&0.3473&0.1982&0.6944&-2.9677&0.9064&0.6714 \\
  0&0&0&0&0&0&0.0921&0.6723&0.2568&0.2794&-3.1073&0.8372 \\
  0&0&0&0&0&0&0.1478&0.4315&0.0098&0.9462&0.0249&-2.5285
\end{array}\!\right] \\
\hline\nonumber
\end{align}
\end{figure*}


\appendix

\subsection{Note on partial minimization}
\label{note_on_parmin}

\noindent The partial minimization is valid for any nonconvex program, i.e.  given any function $F(x,y)$ (not necessarily convex) we always have:
\[   \min_{x,y} F(x,y)  = \min_{x} (\min_y F(x,y)).  \]
Let us briefly prove this statement.  For this, take $x$ fixed and determine, as a function of $x$, that $y$ which minimizes $F(x,y)$ in the second variable. Let us denote this solution  by $y^*(x)  \in  \arg \min_y F(x, y) $ for any fixed $x$. Further, let us define the partial function $f(x) = F(x, y^*(x))$ defined only in the variable $x$. Let $x^* \in \arg \min_x f(x)$. Then, the pair $(x^*, y^*(x^*))$ is an optimal solution of the original problem since for all $(x,y)$ we have:
 \[  F(x,y) \geq F(x,y^*(x))  = f(x) \geq f(x^*) = F(x^*, y^*(x^*)).   \] 
 
\subsection{Parameters of model from Section \ref{sect_ex_net_rand}}
\label{sect_ap_par}

\noindent For the positive network example from Section \ref{sect_ex_net_rand}  the numerical values of the $A$, $B$, $C$ matrices generated randomly in the interval $(-5, 1)$ such that $A$ is a stable matrix and satisfying the interconnection map from Figure  \ref{fig:net_systems} are given in \eqref{eq_ABC_ex_net_rand} and \eqref{eq:BC}. Furthermore, in our numerical experiments we  fix $L=[0\ \ 0 \ \ 0\ \ 1]$.
\begin{equation}\label{eq:BC}
B = \begin{bmatrix}
 0.0569\\0.4503\\0.5825\\0.6866\\0.7194\\0.6500\\0.7269\\0.3738\\0.5816\\0.1161\\0.0577\\0.9798\end{bmatrix},\quad 
C^T=\begin{bmatrix}
 0.2848\\0.5950\\0.9622\\0.1858\\0.1930\\0.3416\\0.9329\\0.3907\\0.2732\\0.1519\\0.3971\\0.3747\end{bmatrix}. \\
\end{equation}


\begin{thebibliography}{100}
\bibitem{antoulas-2005}
A.~C. Antoulas.
\newblock Approximation of large-scale dynamical systems.
\newblock {\em SIAM, Philadelphia}, 2005.

\bibitem{astolfi-TAC2010}
A.~Astolfi.
\newblock Model reduction by moment matching for linear and nonlinear systems.
\newblock {\em IEEE Transactions on Automatic Control}, 50(10): 2321--2336, 2010.

\bibitem{bai-ANM2002}
Z. Bai, 
\newblock Krylov subspace techniques for reduced-order modeling of large-scale dynamical systems.
\newblock {\em Applied Numerical Mathematics}, 43:9--44, 2002.

\bibitem{besselink-sandberg-johansson-ECC2014}
B. Besselink, H. Sandberg,  and K.H. Johansson.
\newblock Model reduction of networked passive systems through clustering.
\newblock In \emph{Proceedings of European Control Conference}, 1069--1074, 2014.

\bibitem{CheKaw:16a}
X. Cheng, Y. Kawano, and J. M. Scherpen. 
\newblock Model reduction of multi-agent systems using dissimilarity-based clustering. 
\newblock  {\em IEEE Transactions on Automatic Control} , DOI: 10.1109/TAC.2018.2853578, 2018.

\bibitem{freund-ACM2004}
R.~W. Freund.
\newblock {SPRIM}: {S}tructure-preserving reduced-order interconnect  macromodeling.
\newblock In {\em Proceedings of Conference on Computer-Aided Design}: 80--87, 2004.

\bibitem{gallivan-vandendorpe-vandooren-SIAM2004}
K.~Gallivan, A.~Vandendorpe, and P.~Van Dooren.
\newblock Model reduction of {MIMO} systems via tangential interpolation.
\newblock {\em SIAM Journal on Matrix Analysis Applications}, 26(2): 328--349, 2004.

\bibitem{gugercin-antoulas-beattie-SIAM2008}
S.~Gugercin, A.~C. Antoulas, and C.~A. Beattie. 
\newblock ${H}_2$ model reduction for large-scale dynamical systems. 
\newblock {\em SIAM Journal on Matrix Analysis \& Applications}, 30:609--638, 2008.

\bibitem{BruRan:14}
C. Grussler and A. Rantzer.
\newblock Modified balanced truncation preserving ellipsoidal cone-invariance.
\newblock In {\em Proceedings of Conference on Decision and Control}, 2014.

\bibitem{IshKas:14} 
T. Ishizaki, K. Kashima, J.  Imura, and K. Aihara. 
\newblock Model reduction and clusterization of large-scale bidirectional networks. 
\newblock {\em IEEE Transactions on Automatic Control}, 59(1): 48--63, 2014.

\bibitem{i-astolfi-colaneri-SCL2014}
T.~C. Ionescu, A.~Astolfi and P.~Colaneri.
\newblock Families of moment matching based low order approximations for linear systems.
\newblock {\em Systems \& Control Letters}, 64: 47--56, 2014.

\bibitem{li-bai-CMS2005}
R.-C. Li and Z. Bai, 
\newblock Structure-preserving model reduction using a Krylov subspace projection formulation. 
\newblock {\em Communications in Mathematical Sciences}, 3:179--199, 2005.

\bibitem{lutowska-2012}
A.~Lutowska.
\newblock  Model order reduction for coupled systems using low-rank approximations.
\newblock {\em PhD thesis}, T.U. Eindhoven,  2012.

\bibitem{MarFra:18}
N. Martin, P. Frasca, and C. Canudas-de-Wit.
\newblock A network reduction method inducing scale-free degree distribution.
\newblock In {\em Proceedings of European Control Conference}: 2236--2241, 2018.

\bibitem{MonTre:13}
N. Monshizadeh, H. L. Trentelman, and M. K. Camlibel. 
\newblock Stability and synchronization preserving model reduction of multi-agent systems.
\newblock{Systems \& Control Letters}, 62(1): 1--10, 2013.

\bibitem{NecNed:11}
I. Necoara, V. Nedelcu, and I. Dumitrache. 
\newblock Parallel and distributed optimization methods for estimation and control in networks. \newblock {\em Journal of Process Control},  21(5): 756--766, 2011.

\bibitem{NecIon:18}
I. Necoara and T. Ionescu. 
\newblock Optimal $H_2$ moment matching-based model reduction for linear systems by (non)convex optimization.
\newblock Tech. Rep., UPB, https://arxiv.org/abs/1811.07409, 2018.

\bibitem{Nes:04}
Yu. Nesterov. 
\newblock  Introductory Lectures on Convex Optimization:  A Basic Course.  
\newblock {\em Kluwer}, 2004.

\bibitem{reis-stykel-MCDMS2007}
T.~Reis and T.~Stykel.
\newblock Stability analysis and model order reduction for coupled systems.
\newblock {\em Mathematical and Computer Modelling of Dynamical Systems}, 13(5): 413--436, 2007.

\bibitem{reis-stykel-SPRINGER2008}
T.~Reis and T.~Stykel.
\newblock A survey on model reduction of coupled systems.
\newblock in W.H.A. Schilders, H.A. van der Vorst, and J. Rommes, editors, {\em Model Order Reduction: Theory, Research Aspects and Applications, Mathematics in Industry,} Springer-Verlag Berlin Heidelberg, vol. 13, pp. 133--155, 2008.

\bibitem{ReiVir:09}
T. Reis and E. Virnik.
\newblock Positivity preserving model reduction.
\newblock in \emph{Positive Systems}, eds. R. Bru and S. Romero-Vivo, LNCIS 389, pp.131--139, Springer, 2009.

\bibitem{sandberg-TAC2010}
H.~Sandberg.
\newblock An extension to balanced truncation with application to structured
  model reduction.
\newblock {\em IEEE Transactions on Automatic Control}, 55(4): 1038--1043, 2010.

\bibitem{sandberg-murray-OCAM2009}
H.~Sandberg and R.~M. Murray.
\newblock Model reduction of interconnected linear systems.
\newblock {\em Optimal Control, Applications and Methods}, 30(3): 225--245, 2009.

\bibitem{SooAnd:14}
A. Sootla and J. Anderson.
\newblock On Projection-Based Model Reduction of Biochemical Networks Part I: The Deterministic Case.
\newblock {\em Conference on Decision and Control}, pp.3615--3620, 2014.

\bibitem{Toi:85}
H.~Toivonen
\newblock A globally convergent algorithm for the optimal constant output feedback problem. 
\newblock {\em International Journal of Control}, 41:1589--1599, 1985.

\bibitem{vandendorpe-vandooren-2009}
A.~Vandendorpe, P.~Van Dooren.
\newblock Model reduction of interconnected systems.
\newblock In {\em Model order reduction: theory, research aspects and applications}: 305--321,  2008.

\bibitem{Ven:06}
A.N.  Venkat.
\newblock Distributed Model Predictive Control: Theory and Applications.
\newblock {\em PhD Thesis},  2006. 
\end{thebibliography}
\end{document}